\documentclass[twoside,a4paper,10pt]{article}
\usepackage{latexsym,amssymb}
\usepackage{mltex}
\usepackage{hyperref}
\usepackage{graphicx}
\usepackage[frenchb,english]{babel}
\usepackage[matrix,arrow]{xy}
\def \be{\begin{eqnarray*}}
\def \ee{\end{eqnarray*}}
\def \ben{\begin{enumerate}}
\def \een{\end{enumerate}}
\def \beit{\begin{itemize}}
\def \eeit{\end{itemize}}
\def \bui#1#2{\mathrel{\mathop{\kern 0pt#1}\limits^{#2}}}
\def \buil#1#2{\mathrel{\mathop{\kern 0pt#1}\limits_{#2}}}
\def \bfll{\begin{flushleft}}
\def \efll{\end{flushleft}}
\def \bflr{\begin{flushright}}
\def \eflr{\end{flushright}}

\def \findemo{\hfill$\square$}
\def \lra{\longrightarrow}
\def \lmt{\longmapsto}
\def \ovl{\overline}

\def \wit{\widetilde}

\def \R{\mathbb{R}}

\def \S{\mathbb{S}}

\newcommand{\rquot}[2]{\raisebox{0.5ex}{$#1$}\!/\!\raisebox{-0.5ex}{$#2$}}

\newtheorem{ethm}{Theorem}[section]

\newtheorem{elemme}[ethm]{Lemma}

\newtheorem{ecor}[ethm]{Corollary}

\title{Almost harmonic spinors}
\author{{\small{Nicolas Ginoux}}\footnote{NWF I - Mathematik,
Universit\"at Regensburg,
D-93040 Regensburg,
\texttt{E-mail: nicolas.ginoux@mathematik.uni-regensburg.de}}, {\small{Jean-Fran\c cois Grosjean}}\footnote{Institut \'Elie Cartan (Math\'ematiques), Universit\'e Henri Poincar\'e Nancy I, B.P.239 F-54506 Vandoeuvre-L\`es-Nancy Cedex,
\texttt{E-mail: grosjean@iecn.u-nancy.fr}}}

\begin{document}
\maketitle

\noindent\begin{center}\begin{tabular}{p{155mm}}
\begin{small}{\bf Abstract.} We show that any closed spin manifold not diffeomorphic to the two-sphere admits a sequence of volume-one-Riemannian metrics for which the smallest non-zero Dirac eigenvalue tends to zero. As an application, we compare the Dirac spectrum with the conformal volume.
\end{small}\\
\end{tabular}\end{center}

\noindent\begin{center}
{\Large Spineurs presque harmoniques}
\vspace{5mm}

\begin{tabular}{p{155mm}}
\begin{small}{\bf R\'esum\'e.} Nous montrons que, sur toute vari\'et\'e spinorielle compacte sans bord non diff\'eomorphe \`a la sph\`ere de dimension deux, il existe une suite de m\'etriques riemanniennes de volume un pour laquelle la plus petite valeur propre non nulle de l'op\'erateur de Dirac tend vers z\'ero. Comme application, nous comparons le spectre de l'op\'erateur de Dirac avec le volume conforme.
\end{small}\\
\end{tabular}\end{center}

\section{Introduction and statement}

\noindent Let $M^n$ be an $n(\geq 2)$-dimensional closed spin manifold and denote by $D_g$ the spin Dirac operator associated to a Riemannian metric $g$. We denote by $\lambda_1(D_g^2)$ and $\lambda_1^+(D_g^2)$ the smallest and the smallest positive eigenvalue of $D_g^2$ respectively. It is well-known that the product $\lambda_1^+(D_g^2)\mathrm{Vol}(M^n,g)^{\frac{2}{n}}$ is scaling-invariant and bounded from below by a positive constant in any conformal class \cite[Thm. 2.3]{Ammspinconf}. One can ask whether the infimum of $\lambda_1^+(D_g^2)\mathrm{Vol}(M,g)^{\frac{2}{n}}$ on the space of all Riemannian metrics remains positive. This holds true if $M$ is the $2$-sphere $\S^2$ since it has only one conformal class; alternatively, it follows from C. B\"ar's estimate \cite{Baersurf} valid for any Riemannian metric $g$ on $\S^2$:
\begin{equation}\label{ineqBaersurf}
\lambda_1(D_g^2)\mathrm{Area}(\S^2,g)\geq4\pi.
\end{equation}
In this respect $\S^2$ is the only exception:

\begin{ethm}\label{texistseq}
For any $n(\geq 2)$-dimensional closed spin manifold $M^n$ not diffeomorphic to $\S^2$ there exists a sequence $(g_p)_{p\in\mathbb{N}}$ of Riemannian metrics on $M^n$ such that $\lambda_1^+(D_{g_p}^2)\mathrm{Vol}(M^n,g_p)^{\frac{2}{n}}\buil{\lra}{p\to\infty}0$.
\end{ethm}

\noindent Therefore one can get the Dirac spectrum as close to $0$ as one wants with fixed volume. Note however that Theorem \ref{texistseq} does not prove the existence of non-zero harmonic spinors, i.e., that $0$ is a Dirac eigenvalue.\\

\noindent Theorem \ref{texistseq} is proved in Section \ref{sproof}. In Section \ref{sappli} we apply it to compare the Dirac spectrum with the conformal volume.\\

\noindent{\bf Acknowledgment.} This work started during a stay of the first-named author at the Institut \'Elie Cartan de Nancy, which he would like to thank for its hospitality and support. It is also a pleasure to thank Bernd Ammann and Emmanuel Humbert for fruitful discussions and their critical reading of the paper.

\section{Proof}\label{sproof}

\noindent The proof of Theorem \ref{texistseq} relies on a standard technique first used in the spinorial context by C. B\"ar \cite{Baer96} to show the existence of metrics with harmonic spinors. Namely we prove the result by gluing a model manifold admitting such a sequence and by studying the convergence of the spectrum on the connected sum. Thus the proof is two-step.

\begin{elemme}\label{lexistseqSnT2}
\noindent\beit\item[i)] {\rm Theorem \ref{texistseq}} holds true on the standard sphere $\S^n$ for any $n\geq 3$.
\item[ii)] {\rm Theorem \ref{texistseq}} holds true on the $2$-torus $\mathbb{T}^2$ endowed with any of its $4$ spin structures.
\eeit 
\end{elemme}

\noindent{\it Proof}:  Both statements follow from elementary arguments.\\
$i)$ For any $n\geq 3$ there exists on $\S^n$ a metric $\wit{g}$ with $\mathrm{Ker}(D_{\wit{g}})\neq \{0\}$: for $n\equiv 3\;(4)$ it is some Berger metric as shown by C. B\"ar \cite[Cor. p.906]{Baer96}, for $n=2m\geq 4$ and $n\geq 5$ the existence of such a metric has been proved by L. Seeger \cite{Seeger00} and M. Dahl \cite[Cor. 4.2]{Dahl06} respectively. Linear interpolation between the standard round metric and $\wit{g}$ provides a smooth one-parameter-family of Riemannian metrics $(g_t)_{0\leq t\leq 1}$ with $g_0=\mathrm{can}$ and $g_1=\wit{g}$. Since the volume remains bounded and the Dirac spectrum depends continuously on the metric in the $\mathrm{C}^1$-topology, we obtain the result.\\
$ii)$ For a real parameter $a>1$ consider the $2$-torus $\mathbb{T}^2:=\rquot{\R^2}{\Gamma}$, where $\Gamma:=\mathbb{Z}\cdot\left(\begin{array}{c}1\\0\end{array}\right)\oplus\mathbb{Z}\cdot\left(\begin{array}{c}0\\a\end{array}\right)$, with induced flat metric $g_a$. It carries $4$ spin structures, $3$ of which can be deduced from each other by an orientation-preserving diffeomorphism of $\mathbb{T}^2$ (see e.g. \cite{AmmHumb06}). Thus it suffices to prove the statement for two spin structures which cannot be obtained from each other by a diffeomorphism, for example for the spin structure inducing a trivial covering on both factors and for the spin structure inducing a trivial covering on the first factor and a non-trivial one on the second one. For the former spin structure the smallest positive eigenvalue of $D_{g_a}^2$ is $\frac{4\pi^2}{a^2}$ and for the latter one it is $\frac{\pi^2}{a^2}$. Since $\mathrm{Area}(\mathbb{T}^2,g_a)=a$ we conclude that in both situations $\lambda_1^+(D_{g_a}^2)\mathrm{Area}(\mathbb{T}^2,g_a)\buil{\lra}{a\to\infty}0$ (compare with \cite[Sec. 3]{AmmHumb06}).\findemo

$ $\\
\noindent In the second step we consider the dimensions $n=2$ and $n\geq 3$ separately. In the latter case it only remains to know how the Dirac spectrum behaves under connected sum\footnote{Thanks to Christian B\"ar for indicating to us the right reference.}:

\begin{ethm}[C. B\"ar \cite{Baer97}]\label{tsurgerycodim3}
Let $(N_1^n,g_1)$ and $(N_2^n,g_2)$ be closed Riemannian spin manifolds of dimension $n\geq 3$.
Let  $L>0$ and $\eta\geq 0$ with $\pm (L+\eta)\notin(\mathrm{Spec}(D_{g_1})\cup\mathrm{Spec}(D_{g_2}))$.\\
Then for any $\varepsilon>0$, there exists a Riemannian metric $\tilde{g}$ on the connected sum $\wit{N}^n:=N_1^n\sharp\, N_2^n$ such that the Dirac eigenvalues of $N_1^n\bui{\cup}{\cdot}N_2^n$ and $(\wit{N}^n,\tilde{g})$ in $]-L-\eta,L+\eta[$ dif\mbox{}fer at most by $\varepsilon$ and that $\mathrm{Vol}(\wit{N}^n,\wit{g})\leq\mathrm{Vol}(N_1^n,g_1)+\mathrm{Vol}(N_2^n,g_2)+\varepsilon$.
\end{ethm}

\noindent Note that, as an easy consequence, Theorem \ref{tsurgerycodim3} remains valid when replacing the eigenvalues of the Dirac operator by those of its square.
Fix now any Riemannian metric $g$ on $M^n$ (with $n\geq 3$). If $p$ is any positive integer, pick from Lemma \ref{lexistseqSnT2} a Riemannian metric $g_p$ of volume one on $\S^n$ with $\lambda_1^+(D_{g_p}^2)\leq\frac{1}{p}$.
Setting $L:=\lambda_1^+(D_{g_p}^2)$, $\varepsilon:=\frac{\lambda_1^+(D_{g_p}^2)}{2}$ and choosing $\eta>0$ with $L+\eta\notin(\mathrm{Spec}(D_{g_p}^2)\cup\mathrm{Spec}(D_g^2))$, Theorem \ref{tsurgerycodim3} implies the existence of a Riemannian metric $\wit{g}_p$ on $\wit{N}^n:=M^n\sharp\,\S^n$ such that at least one eigenvalue of $D_{\wit{g}_p}^2$ lies in the interval $[\frac{\lambda_1^+(D_{g_p}^2)}{2},\frac{3\lambda_1^+(D_{g_p}^2)}{2}]$ and that $\mathrm{Vol}(\wit{N}^n,\wit{g}_p)\leq\mathrm{Vol}(M^n,g)+1+\frac{1}{2p}$. Since $\wit{N}^n$ is spin diffeomorphic to $M^n$ we conclude the proof of Theorem \ref{texistseq} for $n\geq 3$.\\

\noindent In dimension $n=2$ we perform an induction on the genus of the surface. On $\mathbb{T}^2$ Theorem \ref{texistseq} already holds true by Lemma \ref{lexistseqSnT2}. Assume it to hold true for any closed oriented surface $M^2(\gamma)$ of genus $\gamma>0$ and consider a closed oriented surface $M^2(\gamma+1)$ of genus $\gamma+1$. The oriented surface can be obtained as the connected sum of some $M^2(\gamma)$ and $\mathbb{T}^2$. 
Moreover, the spin structure induced on a circle bounding a compact oriented surface is always a non-trivial covering \cite[p.91]{LM}, in particular it itself bounds a disk.
Therefore, every spin structure on $M^2(\gamma+1)$ is induced by some spin structure on $M^2(\gamma)$ and some on $\mathbb{T}^2$\footnote{The argument given in the published version was wrong. Thanks also to Christian B\"ar.}.
It would remain to prove the analog of Theorem \ref{tsurgerycodim3} for surgeries of codimension $2$, at least for connected sums of surfaces. We conjecture this holds true, using arguments and techniques from \cite{AmmDahlHumb06}. Actually much less is needed here:

\begin{elemme}\label{lconvergdim2}
Let $(M_1,g_1)$ and $(M_2,g_2)$ be any oriented closed Riemannian surfaces. Fix $L>0$ and $\eta\geq 0$ with $\pm (L+\eta)\notin(\mathrm{Spec}(D_{g_1})\cup\mathrm{Spec}(D_{g_2}))$.\\
Then for any $\varepsilon>0$ there exists a Riemannian metric $\tilde{g}$ on $M_1\sharp\, M_2$ such that, for any eigenvalue $\lambda$ of $D_{g_1}$ or $D_{g_2}$ in $]-L-\eta,L+\eta[$, there exists an eigenvalue $\tilde{\lambda}$ of $D_{\tilde{g}}$ such that $|\tilde{\lambda}-\lambda|\leq\varepsilon$ and $\mathrm{Area}(M_1\sharp\, M_2,\tilde{g})\leq\mathrm{Area}(M_1,g_1)+\mathrm{Area}(M_2,g_2)+\varepsilon$. 
\end{elemme}

\noindent{\it Proof}: The proof relies on a classical cut-off procedure for eigenvectors of $D_{g_1}$ and $D_{g_2}$. We want to show that $\mathrm{dim}(\mathrm{Ker}(D_{g_1}-\lambda\mathrm{Id}))+ \mathrm{dim}(\mathrm{Ker}(D_{g_2}-\lambda\mathrm{Id}))\leq\mathrm{dim}\Big(\buil{\oplus}{\mu\in[\lambda-\varepsilon,\lambda+\varepsilon]}\mathrm{Ker}(D_{\wit{g}}-\mu\mathrm{Id})\Big)$. Fix $p_i\in M_i$, $i=1,2$, and some sufficiently small $\delta>0$. Consider the connected sum $\wit{M}:=M_1\sharp\, M_2$ obtained by gluing $M_1\setminus B_{p_1}(\delta)$ and $M_2\setminus B_{p_2}(\delta)$ along their boundary, where $B_p(r)$ denotes the open metric disc of center $p$ and radius $r$. From \cite[p.932]{Baer96} or \cite[Sec. 3.1-3.2]{AmmDahlHumb06} there exists a smooth Riemannian metric $\tilde{g}_\delta$ on $\wit{M}$ which coincides with $g_i$ on $M_i\setminus B_{p_i}(\sqrt{\delta})$ and such that $\mathrm{Area}(M_1\sharp\, M_2,\tilde{g}_\delta)\leq\mathrm{Area}(M_1,g_1)+\mathrm{Area}(M_2,g_2)+c\cdot\sqrt{\delta}$, where $c>0$ is a constant depending only on the metrics $g_1$ and $g_2$. In particular we may choose $\delta>0$ such that $c\cdot\sqrt{\delta}<\varepsilon$.
For $i=1,2$ define $\chi_i\in C(\wit{M},[0,1])$ by $\chi_i{}_{|_{M_i\setminus B_{p_i}(\sqrt{\delta})}}:=1$, $\chi_i{}_{|_{B_{p_i}(\delta)}}:=0$ and $\chi_i(x):=2-2\frac{\mathrm{ln}(d(x,p_i))}{\mathrm{ln}(\delta)}$ otherwise, where $d(x,p)$ denotes the distance between $x$ and $p$. Note that $\chi_1$ and $\chi_2$ are well-defined and continuous on the whole $\wit{M}$ and that they can be smoothed out at both $\partial B_{p_i}(\delta)$ and $\partial B_{p_i}(\sqrt{\delta})$ such that the $\mathrm{L}^2$-norm of their gradient changes arbitrarily little. We keep denoting the corresponding smooth functions by $\chi_1$ and $\chi_2$. Consider the map
\be 
\Phi:\mathrm{Ker}(D_{g_1}-\lambda\mathrm{Id})\oplus\mathrm{Ker}(D_{g_2}-\lambda\mathrm{Id})&\lra&\Gamma(\Sigma\wit{M})\\
(\varphi_1,\varphi_2)&\lmt&\chi_1\varphi_1+\chi_2\varphi_2,
\ee
which is well-defined because of $\chi_i{}_{|_{B_{p_i}(\delta)}}=0$ and injective by the unique continuation property (each Dirac eigenvector vanishing on an open subset of a connected Riemannian spin manifold must vanish identically). Now from the min-max principle it suffices to show that $\frac{\|(D_{\tilde{g}_\delta}-\lambda)\varphi\|_{\mathrm{L}^2}}{\|\varphi\|_{\mathrm{L}^2}}\leq\varepsilon$ for all $\varphi\in\mathrm{Im}(\Phi)\setminus\{0\}$. Since the subspaces $\Phi(\mathrm{Ker}(D_{g_1}-\lambda\mathrm{Id}))$ and $\Phi(\mathrm{Ker}(D_{g_2}-\lambda\mathrm{Id}))$ are $\mathrm{L}^2$-orthogonal to each other (for $\mathrm{supp}(\chi_1)\cap\mathrm{supp}(\chi_2)$ has zero measure), we can assume that $\varphi\in\Phi(\mathrm{Ker}(D_{g_1}-\lambda\mathrm{Id}))$ with $\|\varphi\|_{\mathrm{L}^2}=1$. Using the formula $D_g(f\varphi)=df\cdot\varphi+f D_g\varphi$, we compute:
\be 
\int_{\wit{M}}|(D_{\tilde{g}_\delta}-\lambda)\varphi|^2 v_{\tilde{g}_\delta}&=&\int_{M_1}|(D_{\tilde{g}}-\lambda)\chi_1\varphi_1|^2v_{\tilde{g}_\delta}\\
&=&\int_{M_1}|d\chi_1|^2|\varphi_1|^2v_{\tilde{g}_\delta}\\
&\leq&C \sup_{M_1}(|\varphi_1|^2)\int_\delta^{\sqrt{\delta}}\frac{4}{r^2\mathrm{ln}(\delta)^2}rdr\\
&\leq&-\frac{C'}{\mathrm{ln}(\delta)},
\ee
where $C>0$ is a constant depending only on the original metrics in the ring $B_{p_1}(\sqrt{\delta})\setminus B_{p_1}(\delta)$ and $C'=2C\cdot\buil{\sup}{\varphi_1\in\mathrm{Ker}(D_{g_1}-\lambda\mathrm{Id})}(\buil{\sup}{M_1}(|\varphi_1|^2)$ (note that $C'<\infty$ since $\mathrm{Ker}(D_{g_1}-\lambda\mathrm{Id})$ is finite-dimensional). We deduce that $\|(D_{\tilde{g}_\delta}-\lambda)\varphi\|_{\mathrm{L}^2}^2\buil{\lra}{\delta\to 0}0$ and the statement of Lemma \ref{lconvergdim2}.
\findemo

$ $\\

\noindent The proof of Theorem \ref{texistseq} for $n=2$ follows the lines of that for $n\geq 3$: given any Riemannian metric $g$ on $M^2(\gamma)$ and a positive $p\in\mathbb{N}$, pick from Lemma \ref{lexistseqSnT2} a Riemannian metric $g_p$ of unit area on $\mathbb{T}^2$ with $\lambda_1^+(D_{g_p}^2)\leq\frac{1}{p}$, whatever the spin structure of $\mathbb{T}^2$ is. Lemma \ref{lconvergdim2} ensures the existence of a Riemannian metric $\wit{g}_p$ on $M^2(\gamma+1)=M^2(\gamma)\sharp\,\mathbb{T}^2$ such that at least one eigenvalue of $D_{\wit{g}_p}^2$ lies in the interval $[\frac{\lambda_1^+(D_{g_p}^2)}{2},\frac{3\lambda_1^+(D_{g_p}^2)}{2}]$ and that $\mathrm{Area}(M^2(\gamma+1),\wit{g}_p)\leq\mathrm{Area}(M^2(\gamma),g)+1+\frac{1}{2p}$. This proves the result for $\gamma+1$ and concludes the proof of Theorem \ref{texistseq}.

\section{Application}\label{sappli}

\noindent This note was motivated by the study of the relationship between the Dirac spectrum and the so-called conformal volume, which is the conformal invariant defined for any closed Riemannian manifold $(M^n,g)$ by
\[\mathrm{V}_c(M^n,[g]):=\inf_{N\in\mathbb{N}}\Big(\inf_{\varphi\in\mathrm{Imm}_c(M^n,\S^N)}\big(\sup_{\gamma\in\mathrm{Conf}(\S^N)}(\mathrm{Vol}(M^n,(\gamma\circ\varphi)^*\mathrm{can}))\big)\Big),\]
where $\mathrm{Imm}_c(M^n,\S^N)$ denotes the set of conformal immersions $(M^n,g)\lra(\S^N,\mathrm{can})$ and $\mathrm{Conf}(\S^N)$ the group of conformal diffeomorphisms of $(\S^N,\mathrm{can})$. First introduced by P. Li and S.-T. Yau \cite{LiYau82}, it has been shown to be directly related to the Laplace spectrum since it provides an upper bound of the corresponding spectral invariant \cite{LiYau82,ElSI86}: $(0<)\lambda_1(\Delta)\mathrm{Vol}(M^n,g)^{\frac{2}{n}}\leq n\mathrm{V}_c(M^n,[g])^{\frac{2}{n}}$. For the Dirac operator such a result cannot be expected because of $\buil{\sup}{\ovl{g}\in [g]}\big(\lambda_1^+(D_{\ovl{g}}^2)\mathrm{Vol}(M^n,\ovl{g})^{\frac{2}{n}}\big)=\infty$, see \cite[Thm. 1.1]{AmmJammes07}. However, one could reasonably conjecture that the conformal volume bounds $\lambda_1(D_g^2)\mathrm{Vol}(M,g)^{\frac{2}{n}}$ from below, provided the possible eigenvalue $0$ is left aside. For $M=\S^2$ this is the case because of (\ref{ineqBaersurf}) and $4\pi=\mathrm{Area}_c(\S^2)$ (see \cite{LiYau82}). It is hopeless for any other manifold:

\begin{ecor}\label{cconflbcv}
For any $n(\geq 2)$-dimensional closed Riemannian spin manifold $(M^n,g)$ not diffeomorphic to $\S^2$ there exists no positive constant $c(M)$ (depending only on $M$) such that 
\[\lambda_1^+(D_g^2)\mathrm{Vol}(M,g)^{\frac{2}{n}}\geq c(M)\mathrm{V}_c(M^n,[g])^{\frac{2}{n}}. \]
\end{ecor}

\noindent{\it Proof}: It is elementary to show that \cite[Fact 2]{LiYau82}
\begin{equation}\label{ineqLiYau}
\mathrm{V}_c(M^n,[g])\geq\mathrm{Vol}(\S^n,\mathrm{can}), 
\end{equation}
whose r.h.s. does not depend on the metric $g$. We conclude with Theorem \ref{texistseq}.
\findemo
 $ $\\

\noindent Still there exists a subtle relationship between the Dirac spectrum and the conformal volume. Indeed by \cite[Thm. 3.1 \& 3.2]{Ammspinconf} and \cite[Thm. 1.1]{AmmGrosjHumbMorel}, $\inf_{\ovl{g}\in[g]}\Big(\lambda_1^+(D_{\ovl{g}}^2)\mathrm{Vol}(M^n,\ovl{g})^{\frac{2}{n}}\Big)\leq\frac{n^2}{4}\mathrm{Vol}(\S^n,\mathrm{can})^{\frac{2}{n}}$,
hence combining with (\ref{ineqLiYau}) one obtains
\[\inf_{\ovl{g}\in[g]}\Big(\lambda_1^+(D_{\ovl{g}}^2)\mathrm{Vol}(M^n,\ovl{g})^{\frac{2}{n}}\Big)\leq\frac{n^2}{4} \mathrm{V}_c(M^n,[g])^{\frac{2}{n}}.\]

\providecommand{\bysame}{\leavevmode\hbox to3em{\hrulefill}\thinspace}

\end{document}